\documentclass[11pt]{amsart}

\pdfoutput=1
\usepackage{amssymb,verbatim}
\usepackage[colorlinks=true, pdfstartview=FitV, linkcolor=blue, citecolor=blue, urlcolor=blue]{hyperref} 
\usepackage{graphicx}
\usepackage{subfigure,wrapfig}

\newcommand{\R}{{\mathbb R}}
\newcommand{\N}{{\mathbb N}}

\newtheorem*{thm}{Theorem}\theoremstyle{definition}

\everymath{\displaystyle}

\begin{document}
\title{A fractal version of the pinwheel tiling}
\author{Natalie Priebe Frank}
\address{Department of Mathematics, Vassar College, Box 248, Poughkeepsie, NY  12604, USA}
\email{nafrank@vassar.edu}

\author{Michael F. Whittaker }
\address{School of Mathematics and Applied Statistics,
University of Wollongong, Wollongong, NSW, Australia 2500}
\email{mfwhittaker@gmail.com}
\thanks{The authors would like to thank Dirk Frettl\"oh for helpful discussions about the aorta and Edmund Harriss for pointing out our theorem on rotations in the aorta. The second author is partially supported by ARC grant 228-37-1021, Australia}

\maketitle

\centerline{\em Dedicated to the inspiration of Benoit Mandelbrot.}

\setlength{\baselineskip}{.8cm}
The pinwheel tilings are a remarkable class of tilings of the plane, and our main goal in this paper is to introduce a fractal version of them.  We will begin by describing how to construct the pinwheel tilings themselves and by discussing some of the properties that have generated so much interest.   After that we will develop the fractal version and discuss some of its properties.  Finding this fractile version was an inherently interesting problem, and the solution we found is unusual in the tiling literature.

Like the well-known Penrose tilings \cite{Gardner}, pinwheel tilings are generated by an ``inflate-and-subdivide rule" (see Figure \ref{orig.pin.subs}). Tilings generated by inflate-and-subdivide rules form a class of tilings that have a considerable amount of global structure called self-similarity.  Self-similar tilings are usually nonperiodic but still exhibit  a form of  ``long-range order" that makes their study particularly fruitful.   Unlike Penrose tilings and most known examples of self-similar tilings, tiles in any pinwheel tiling appear in infinitely many different orientations.   The pinwheel tilings were the first example of this sort and as such presented both new challenges and intriguing properties.

Many examples of self-similar tilings are made of {\em fractiles:}  tiles with fractal boundaries.   Fractiles arise in the foundational work \cite{Kenyon} for constructing a self-similar tiling for a given inflation factor.   Two fractile versions of the Penrose tilings are introduced in \cite{Bandt-Gummelt}. Additionally, the procedure used in \cite{Me.Boris} may result in self-similar tilings made up of fractiles.    This made it reasonable to expect that the
pinwheel tilings might have a fractal variant, but did not provide a template for finding it.

The technique for finding fractiles in both \cite{Bandt-Gummelt} and \cite{Me.Boris} is similar. One begins with an inflate-and-subdivide rule for which the edges of each inflated tile do not quite match up with the edges of the tiles that replace it (the  Penrose kite and dart are an example).  The edges of each tile are redrawn using the edges of the tiles that replace it.   These new edges are revised iteratively by the subdivision rule ad infinitum.   The final result is a set of fractiles that are redrawings of the original tiles, but now they inflate-and-subdivide perfectly.    Our technique is completely different.   We found a fractal that runs through the {\em interior} of a pinwheel triangle and behaves nicely under the inflate-and-subdivide rule.     The fact that this fractal extends to become the boundary of fractiles follows naturally (with some work) from the pinwheel inflation, but there was no way to know {\em how many types} of fractiles to expect.   It was only by creating the images by computer that we were able to
generate enough information to answer that question.\footnote{Mathematica code for the images in this paper is available on request by contacting the first author.}


\section*{Pinwheel tilings}
Pinwheel tilings are made up of right triangles of side lengths 1, 2, and $\sqrt{5}$.  We say a pinwheel triangle is in  {\em standard position} and call it a {\em standard triangle} if its vertices are at $(-.5, -.5), (.5, -.5)$, and $(-.5, 1.5)$.   If we multiply this standard triangle by the matrix $M_P=$ {\small $\left( \begin{array}{cc} 2 & 1 \\ -1 & 2 \end{array} \right)$} , it can be subdivided into five pinwheel triangles of the original size (see Figure 1).   This is known as the pinwheel {\em inflate-and-subdivide rule}, or more simply as the pinwheel {\em substitution rule} \cite{Rad}.
(Readers who wish to make drawings for yourselves: notice that all of the images in this paper are oriented with a standard triangle at the origin and the origin marked.)  
\begin{figure}[ht]
\parbox{.5in}{\includegraphics[width=.35in]{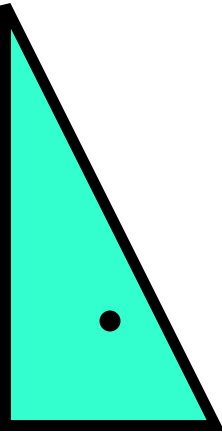}}
\raisebox{-.2in}{\text{\LARGE $\to$}}
\raisebox{.2in}{\parbox{1.2in}{\includegraphics[width=.75in]{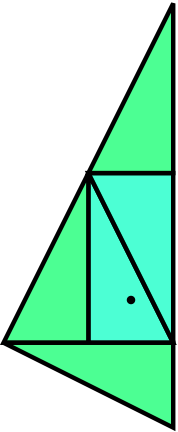}}}
\caption{The pinwheel inflate-and-subdivide rule.}
\label{orig.pin.subs}
\end{figure}
We can apply the rule again, multiplying by $M_P$ and then subdividing each of the five tiles as in Figure \ref{orig.pin.subs} and its reflection.  In this way we obtain a patch of 25 tiles that we call a {\em level-2 tile}; substituting $n$ times produces a {\em level-$n$ tile}.  In Figure \ref{orig.pin.its} we see three levels of the substitution, where we have emphasized the borders of the original five tiles to exhibit the hierarchy.

\begin{figure}[ht]
\centerline{
\raisebox{.05in}{\parbox{.35in}{\includegraphics[width=.35in]{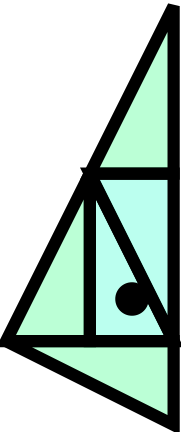}}}
\raisebox{0in}{\text{\LARGE $\to$}}
\raisebox{.1in}{\parbox{1.5in}{\includegraphics[width=1.35in]{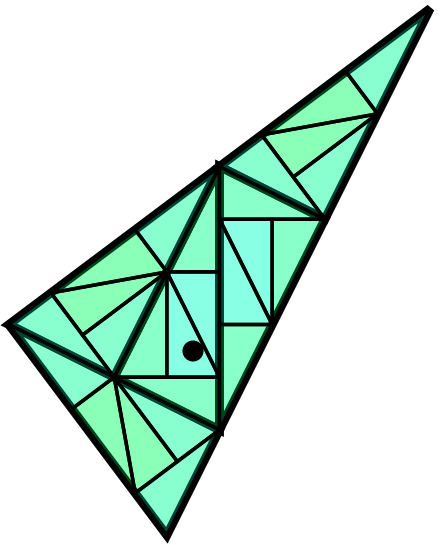}}}
\raisebox{0in}{\text{\LARGE $\to$}}
\parbox{3.5in}{\includegraphics[width=3.6in]{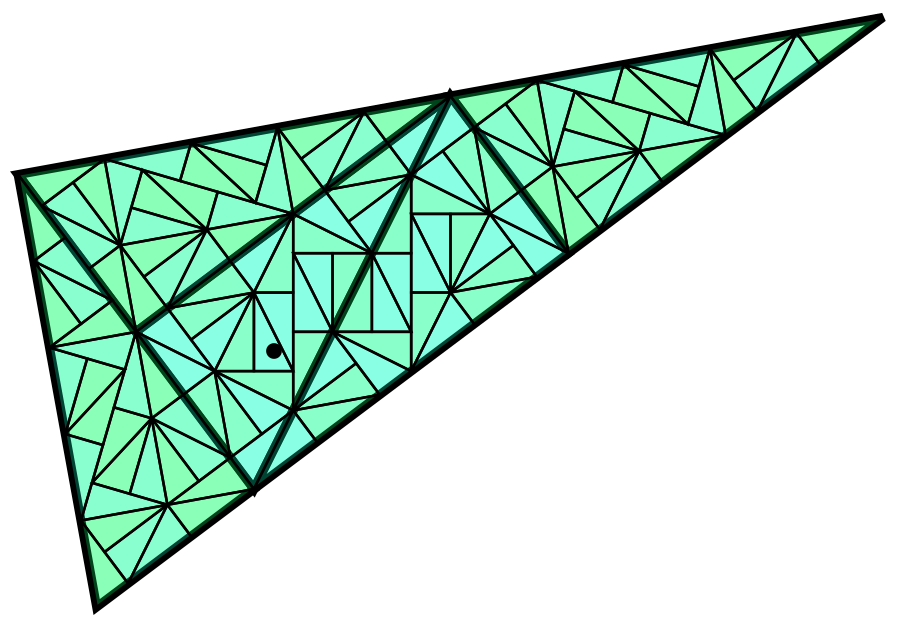}}
}
\caption{Level-1, -2, and -3 tiles for the pinwheel inflate-and-subdivide rule.}
\label{orig.pin.its}
\end{figure}

Many inflate-and-subdivide rules for tilings have been discovered since attention was first drawn to the subject in the 1960's.
A compendium of tilings generated by substitution rules appears on the Tilings Encyclopedia website \cite{Tiling.encyclopedia}, and an introduction to several different forms of tiling
substitutions appears in \cite{My.primer}.

There are several equivalent ways to obtain infinite tilings from a tiling substitution rule
.  The most straightforward is a constructive approach.  Since the standard triangle is invariant under the pinwheel inflation, the level-1 tile it becomes will be invariant under any further applications of the inflate-and-subdivide rule.   So will the level-$n$ tiles once the rule has been iterated at least $n$ times.    Thus it is easy to see that when the substitution rule has been applied ad infinitum we will obtain a well-defined infinite tiling $T_0$ of the plane.    In fact we can apply the inflate-and-subdivide rule to {\em any} tiling of the plane made up of pinwheel triangles by multiplying by the matrix $M_P$ and then subdividing.  The tilings that are invariant are called self-similar tilings and one can show that they are all rotations of $T_0$.   

But there are other infinite tilings we would like to call pinwheel tilings.   If we slide $T_0$ so that the origin is in some other tile, or if we rotate or reflect $T_0$, or if we apply any rigid motion to all the tiles in $T_0$ we have really only changed the placement of $T_0$ in $\R^2$.  Thus we will consider any translation or other rigid motion of $T_0$ to be a pinwheel tiling.  Moreover, if a tiling has the property that the patch of tiles in every large ball around the origin agrees with that in some rigid motion of $T_0$, we will call it a pinwheel tiling as well.  We consider all infinite pinwheel tilings we just described  to be elements in the {\em tiling space} $X_P$.

We can summarize this approach to defining pinwheel tilings as follows:  1) construct $T_0$, a pinwheel tiling that is self-similar;
2) Act on $T_0$ by all possible rigid motions, obtaining infinitely many `different' pinwheel tilings; and 3) include any other tilings that agree with tilings from step 2) over all arbitrarily large but finite regions of the plane.    Every pinwheel tiling in $X_P$ looks locally like $T_0$, but there are infinitely many tilings obtained in step 3) that are not rigid motions of $T_0$.    In fact these tilings are accumulation points of the tilings from step 2) under a ``big ball" metric that says that two tilings are close if they very nearly agree on a big ball around the origin.   Including the tilings in step 3) makes $X_P$ a compact topological space.

As noted earlier, one of the main reasons that pinwheel tilings are of such importance is that in any pinwheel tiling, the triangles appear in a countably infinite number of distinct orientations.   This isn't difficult to see once one notices that the pinwheel angle $\phi = \arctan(1/2)$ is irrational with respect to $\pi$, and governs the orientations we see in Figure \ref{orig.pin.its}.   This leads to the fact that the space of pinwheel tilings can be decomposed into the product of an oriented tiling space and a circle \cite{Whi}.

\subsection*{Two intriguing pinwheel properties}
Like many tiling spaces generated from inflate-and-subdivide rules, the pinwheel space has a sort of homogeneity known as {\em unique ergodicity} \cite{Rad1}.   In the situation where the tiles appear in only finitely many rotations, unique ergodicity automatically implies that every finite configuration of tiles appears with a well-defined frequency in every tiling $T$ in the tiling space.    The frequency of some patch $C$ of tiles can be computed by looking at the number of times $C$ occurs in some large ball in $T$ and dividing by the area of that ball; the fact that there will be a limit as the size of the balls goes to infinity is a result of unique ergodicity.   This approach doesn't quite work when there are infinitely many rotations in every tiling.

Unique ergodicity in the pinwheel case means that there is a statistical form of rotational invariance present in $X_P$ that is quite intriguing.   Consider a finite configuration $C$ of tiles and some interval $I$ of orientations in which it might appear.   Given a tiling $T \in X_P$ we can count the number of times that $C$ appears in a large ball in $T$ in an orientation from $I$.   Dividing that by the size of the ball, and then taking a limit as the size of the balls goes to infinity gives the frequency of occurrence of $C$ in orientation $I$.   The fact that this frequency is independent not only of the tiling chosen but of the sequence of balls in $T$ is a side effect of unique ergodicity.   What is more remarkable, the frequency depends only on the {\em size} of $I$, not on $I$ itself \cite{Radin.miles}.   Thus not only are the rotations uniformly distributed, no particular range of orientations is preferred over another.   For this reason the pinwheel tiling space is considered ``statistically round" even though most individual tilings in it are not rotationally invariant.

Another surprising property of pinwheel tilings is that the hierarchical structure mandated by the inflate-and-subdivide rule can be enforced by local constraints called {\em matching rules} \cite{Rad}, decorations on the edges of tiles that specify how they are allowed to meet up.    Although many famous tilings, for instance the Penrose tilings, were known to come equipped with matching rules that force the hierarchical structure, this was the first example for which the matching rules also enforced infinite rotations.  In \cite{Rad}, a new set of triangles is constructed by making numerous copies of the pinwheel triangles, each with markings on their edges that specify how they are allowed to meet.    The remarkable fact is that this extremely local constraint forces the pinwheel hierarchy: any tiling with these new triangles that obeys the matching rules will become a tiling from $X_P$ when the markings on the edges are forgotten.   

\subsection*{The kite-domino version of pinwheel tilings}
A useful concept in tiling theory is that of {\em mutual local derivability}, which gives a way of comparing tilings built with different tile shapes.  Given two tilings $T_1$ and $T_2$ of $\R^2$,  we say that $T_2$ is {\em locally derivable} from $T_1$ if there is a finite radius $R$ such that the patch in the ball of radius $R$ about any point $\vec{x} \in \R^2$ determines the precise type and placement of the tile (or tiles) in $T_2$ at $\vec{x}$.  If $T_2$ is locally derivable from $T_1$ and $T_1$ is also locally derivable from $T_2$, we say the tilings are mutually locally derivable.   If two tiling spaces are mutually locally derivable, then they are homeomorphic in the big ball topology.  The main goal of this work is to introduce a tiling substitution on fractal tiles that produces tilings that are mutually locally derivable from the pinwheel tilings.   But first, following \cite{BFG2}, we introduce a tiling substitution called the ``kite-domino" pinwheel tilings.

\begin{figure}[h]
\includegraphics[width=.75in]{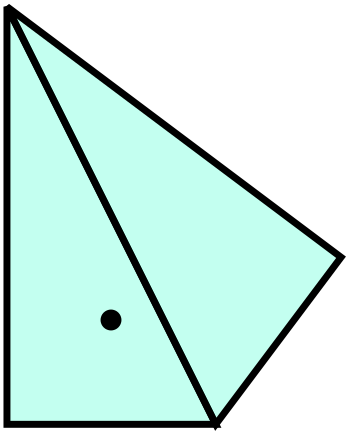} \hspace{.5cm} \includegraphics[width=.5in]{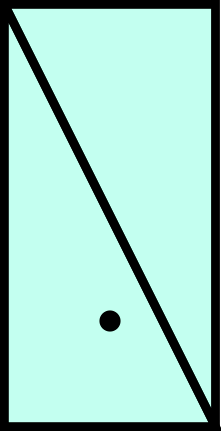}
\caption{The kite and domino}
\label{kite}
\end{figure}

\begin{figure}[h]
\includegraphics[width=2.1in]{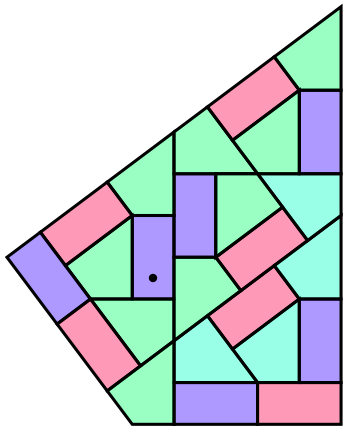} \hspace{1cm} \includegraphics[width=3in]{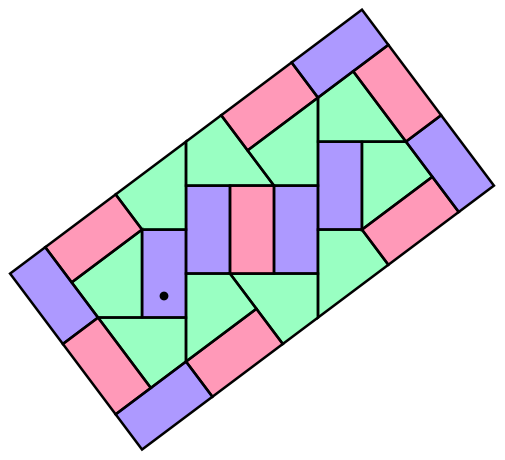}
\caption{Substitution for the kite and domino}
\label{kd.subs}
\end{figure}

The pinwheel triangles in any pinwheel tiling meet up hypotenuse-to-hypotenuse to form either a {\em kite} or a {\em domino}, which we show in standard position in Figure \ref{kite}.   There are two types of domino: the one pictured and one with an opposite diagonal; in our images we denote the difference by shading them differently.  
It is clear that every tiling in $X_P$ can be locally transformed into a kite-domino tiling by fusing together triangles along each hypotenuse.  If the pinwheel substitution is applied to a kite or domino twice the result can be composed into kites and dominoes, resulting in the inflate-and-subdivide rule of Figure \ref{kd.subs}.

We can build the space $X_{KD}$ of all tilings admitted by the kite-domino substitution using the same three-step `constructive' method we used to define the pinwheel tiling space $X_P$.   It is shown in \cite{BFG2} that the pinwheel tiling space $X_P$ is mutually locally derivable from the kite-domino substitution tiling space $X_{KD}$.   In what follows we will rely on $X_{KD}$ to make the fractal version of the pinwheel tilings.

\section*{Construction of pinwheel fractiles}
\label{fractile.sec}

To construct the pinwheel fractiles we construct a fractal, invariant under substitution, that we will use to mark all pinwheel triangles. We call this fractal the {\em aorta}. The aorta will be used both to form the boundaries of the fractiles and to define the local map taking pinwheel tilings to fractal pinwheel tilings.

\subsection*{The aorta}
\label{aorta.sec}
 
\begin{figure}[ht]
\centering
\subfigure{\includegraphics[width=.75in]{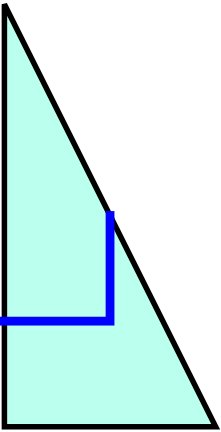}}
\raisebox{.5in}{\text{\LARGE $\to$}}
\subfigure{\includegraphics[width=.75in]{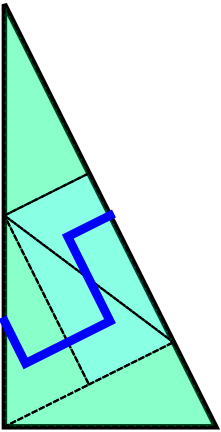}}
\raisebox{.5in}{\text{\LARGE $\to$}}
\subfigure{\includegraphics[width=.75in]{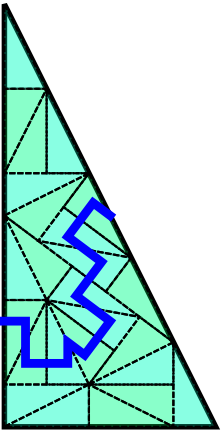}}
\raisebox{.5in}{\text{\LARGE $\to$}}
\subfigure{\includegraphics[width=.75in]{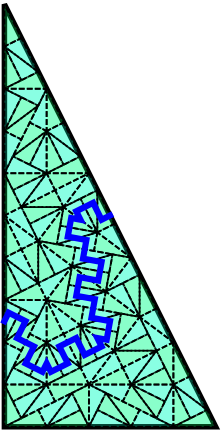}}
\raisebox{.5in}{\text{\LARGE $ \cdots \to$}}
\subfigure{\includegraphics[width=.75in]{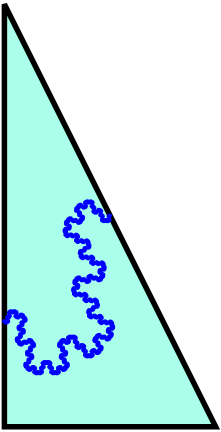}}

\caption{The subdivision method for generating the aorta.}
\label{subdiv.to.fractal}
\end{figure}

There are three special points in a pinwheel triangle: the origin, the point $(-.5, 0)$, and the point $(0, .5)$.   The origin is a {\em (central) control point} (cf. \cite{Sol.self.similar}) since its location in the triangle is invariant under substitution.   We will call the points $(-.5,0)$ and $(0,.5)$ the {\em side} and {\em hypotenuse control points}, respectively.   The key observation is that one can generate a fractal by connecting these three control points and then iterating the pinwheel subdividision rule without inflating.
Figure \ref{subdiv.to.fractal} shows a sequence of subdivisions of the standard triangle.  The side and hypotenuse control points alternate type in the subtriangles. The resulting fractal is the aorta.

Alternatively, one can define the aorta to be the invariant set of an iterated function system.  Let $M_P$ be the pinwheel expansion matrix given above, and let $R_y$ and $R_\pi$ denote reflection across the $y$-axis and rotation by $\pi$, respectively.
Let $f_1(x,y) = M_P^{-1} * R_y(x,y) + (-0.4,-0.2)$, $f_2(x,y)= M_P^{-1}(x,y)$, and $f_3(x,y) = R_\pi*M_P^{-1}(x,y) + (-0.2,0.4)$.   Note that the union of these maps take each stage of the aorta to the next in Figure \ref{subdiv.to.fractal}. Since each $f_j$ is a contraction, there is a unique set such that $A = \bigcup^3_{j = 1} f_j(A)$, and of course $A$ is the aorta.    
 
\begin{figure}[ht]
\centerline{
\raisebox{.05in}{\parbox{.35in}{\includegraphics[width=.35in]{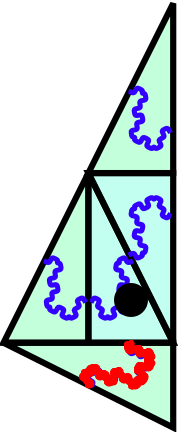}}}
\raisebox{0in}{\text{\LARGE $\to$}}
\raisebox{.1in}{\parbox{1.5in}{\includegraphics[width=1.35in]{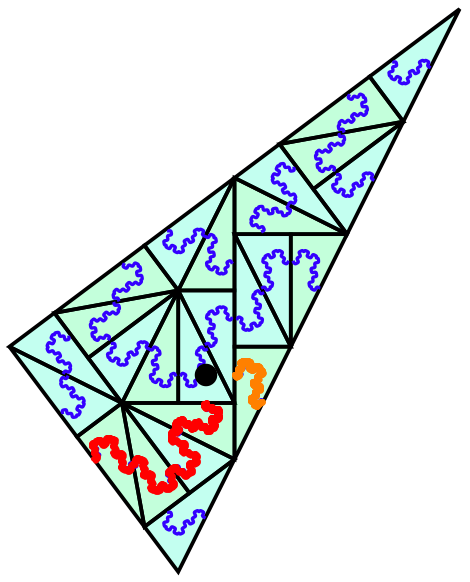}}}
\raisebox{0in}{\text{\LARGE $\to$}}
\parbox{3.5in}{\includegraphics[width=3.6in]{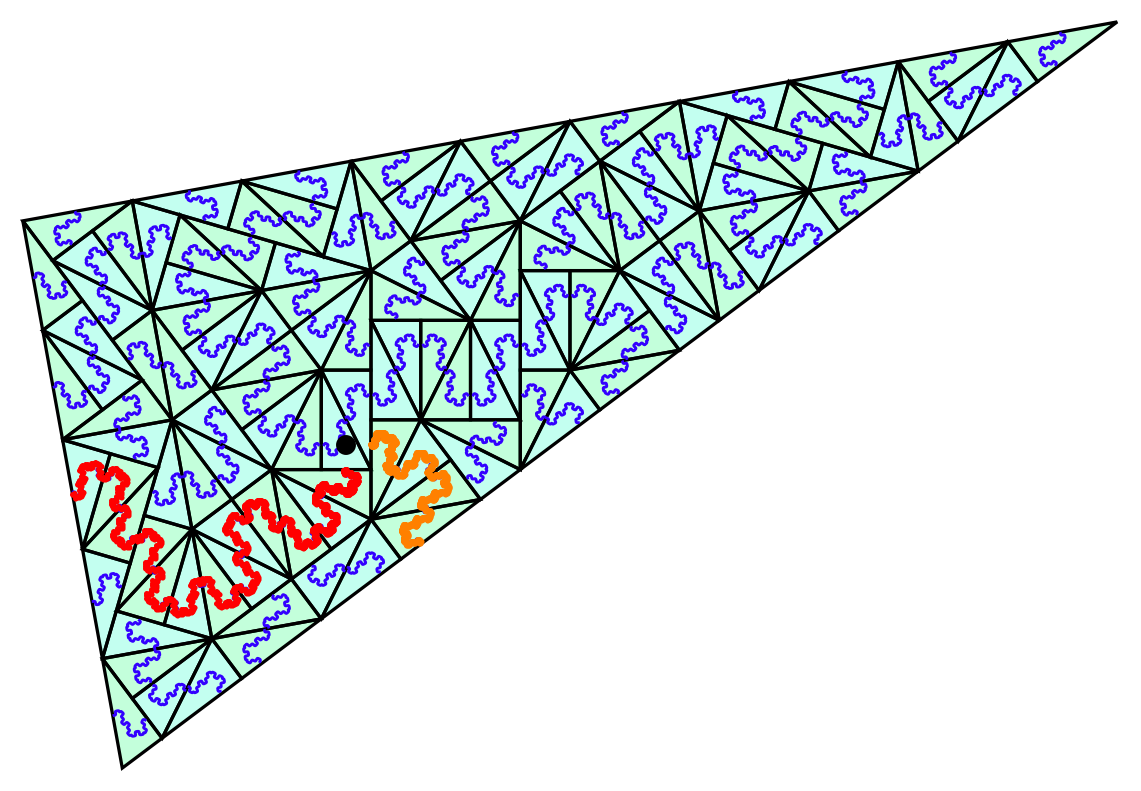}}
}
\caption{A few iterations of the pinwheel inflation with the aortas marked.}
\label{aorta.markings}
\end{figure}

If we begin with a pinwheel triangle in standard position, mark it with its aorta, and inflate by $M_P$, then the aorta will lie along the aortas of three of the five tiles in its subdivision.  So what happens if we mark the aortas of all five tiles in this level-1 tile?   Upon substitution, these five aortas will lie atop fifteen of the 25 aortas in the level-2 tile, all shown in Figure \ref{aorta.markings}.  
A close look at the marked level-3 tile of Figure \ref{aorta.markings} suggests two things.   First, that it may be possible to join up the aortas to create a finite set of tiles with fractal boundary.   Second, the forward invariance of the aorta sets of level-$N$ tiles indicates that these fractiles may possess their own substitution rules.  We will show that both are true.    But first, a few questions and some discussion about the aortas themselves.

\noindent{\sc Questions.} Suppose we mark the aortas of a pinwheel tiling $T$ in $X_P$.   Each individual aorta is part of a connected component of aortas that crosses $N$ triangles, where $N$ is either a positive integer or infinity.   
\begin{enumerate}
\item What is the distribution of $N$ over the tiling $T$?\footnote{We thank one of the reviewers for this very interesting question.}
\item The infinite tiling obtained by continuing to substitute the pinwheel tiling pictured in Figures \ref{orig.pin.its} and \ref{aorta.markings} has a two-sided infinite aorta passing through the origin.
Does every pinwheel tiling have a one- or two-sided  infinite aorta?   If not, what proportion of tilings do?
\end{enumerate}

 \noindent{\sc Note on generalizing the aorta.}  The fact that the side and hypotenuse control points are related by substitution from one step to the next is essential to the existence of the aorta. When we subdivide a pinwheel tile, the hypotenuse control point becomes a side control point and vice versa; an arbitrary finite set of points on the boundary will not behave so nicely.  Finding `aortas' (hidden fractals) in other tiling substitutions using our method will involve finding points on the tile boundaries that are related by substitution.   The set of prototile vertices seems like a good place to start looking, but we have not yet been able to discover any nontrivial examples of these hidden fractals in other tiling substitutions.  It would be interesting to look deeper
into the issues determining which tiling substitutions have versions of the aorta of their own.
 
 \subsection*{The fractiles}
 \label{fractile.sec}

There are two ways to construct equivalent (up to rescaling) versions of the fractal pinwheel tiles.  One way is to begin with the aorta marking of Figure \ref{aorta.markings} and join the aortas that stop abruptly at a tile edge to the central control point of the adjacent tile using an appropriate fractal (a piece of the aorta, in fact).   We call this the {\em continuation method}.  Alternatively, we can mark the pinwheel tiles more elaborately, marking not the aorta, but instead the five aortas of the tiles in the subdivision of each tile.   Connecting the dangling aortas to nearby control points can be done unambiguously with kites and dominoes, and we call this the {\em kite-domino method}.   The tiles produced by the continuation method are equivalent to those from the kite-domino method except that they are five times as large.

\subsubsection*{The continuation method}
In any pinwheel tiling the triangles meet full hypotenuse to full hypotenuse; this implies that whenever an aorta does not connect to an adjacent aorta, it is at the side control point.  The surprising fact is that there are only two ways that this can happen.   In Figure \ref{aorta.markings} we have highlighted (in red and orange) one instance of each type of dangling aorta and follow how each continues after substitution.   
We call the two types of continuations these require the {\em main} and {\em domino} continuations, respectively.    (Notice that the triangular patch for the domino continuation does not appear until the second substitution.)

\begin{figure}[ht]
\centering
\subfigure{\includegraphics[width=1in]{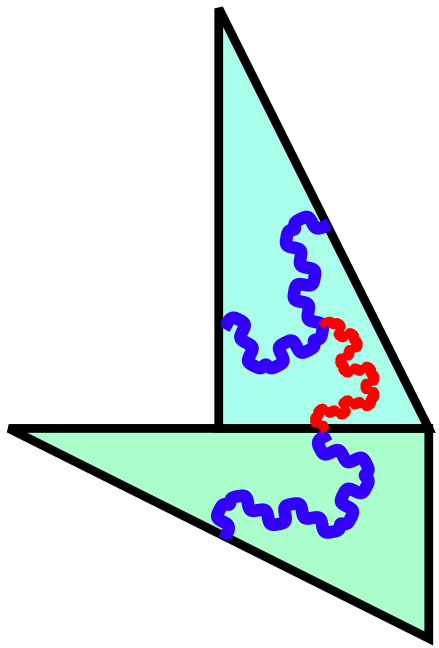}}
\hspace{1cm}
\subfigure{\includegraphics[width=1in]{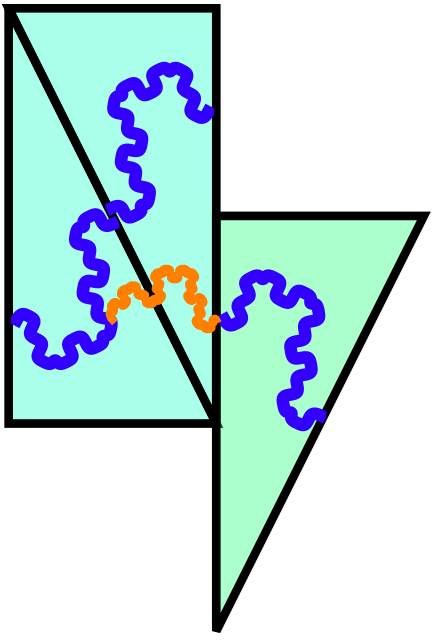}}
\caption{The two types of continuations.}
\label{dangling}
\end{figure}

Further substitution indicates how to
define the continuations: they are isometric copies of the part of the aorta connecting the side control point to the central control point.  The continuations appear in their triangle patches as  pictured in Figure \ref{dangling}, with the main one on the left  and the domino one on the right (red and orange added for emphasis only).  That both kinds of continuation behave well under substitution can be seen from an iterated function systems argument similar to that for the aorta.

Given any pinwheel tiling $T$, we can now produce a new tiling $T_C$ with fractal boundary by marking all aortas as in Figure \ref{aorta.markings} and then adding the continuations as prescribed by Figure \ref{dangling}.   We defer discussion of the properties of tilings produced by this method, preferring to discuss the equivalent tilings produced by the kite-domino method.

\subsubsection*{The kite-domino method}
\begin{figure}[ht]
\centering
\subfigure[]{
\includegraphics[height=1.55in]{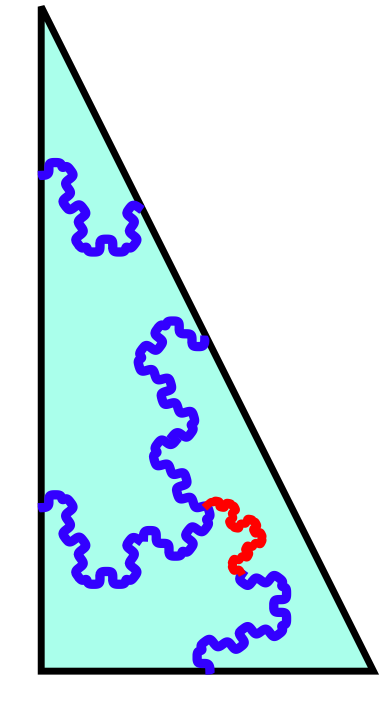}
\label{mark.aorta}
}
\subfigure[]{
\includegraphics[height=1.55in]{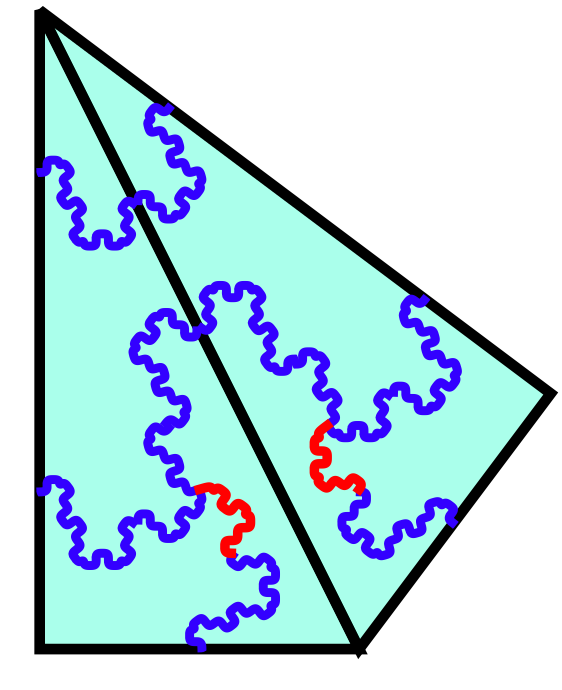}
\label{mark.kite.ventricle}
}
\subfigure[]{
\includegraphics[height=1.55in]{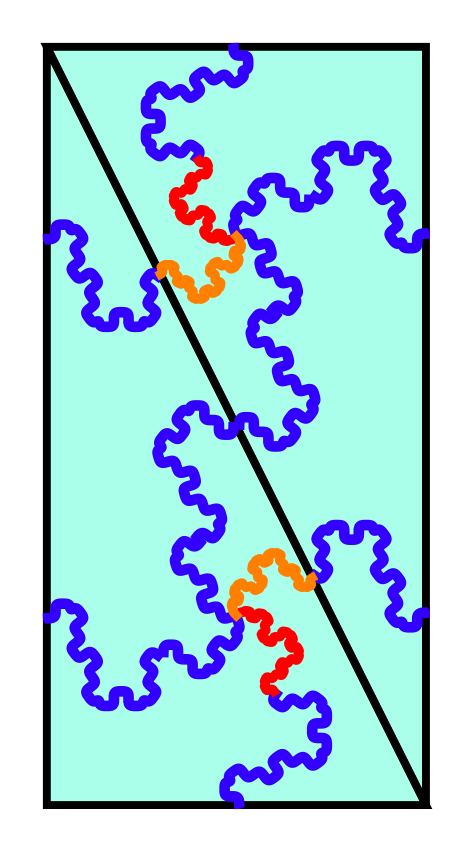}
\label{mark.domino.ventricle}
}
\caption{Marking the kite and domino}
\label{aortas.ventricles.kd}
\end{figure}

We begin with a pinwheel triangle, marking not its aorta but instead the five {\em sub-aortas} that are the preimages of the aortas in its level-1 triangle.  We must add an additional fractal segment to connect the dangling sub-aorta to the central control point.  This segment is shown in red in Figure \ref{mark.aorta}; the dangling sub-aorta along with this segment form exactly the main continuation shown on the left of Figure \ref{dangling}.  The marking of the kite tile, shown in Figure \ref{mark.kite.ventricle}, is simply this initial marking on both of its triangles.   In order to mark a domino tile, we need to use the initial markings on its two triangles, but we also need to resolve the two dangling sub-aortas that arise along the hypotenuse.   As shown in Figure \ref{mark.domino.ventricle}, we add fractal segments to connect these to the central control points so that the resulting fractals are the domino continuations of Figure \ref{dangling}.

\begin{figure}[ht]
\includegraphics[width=2.75in]{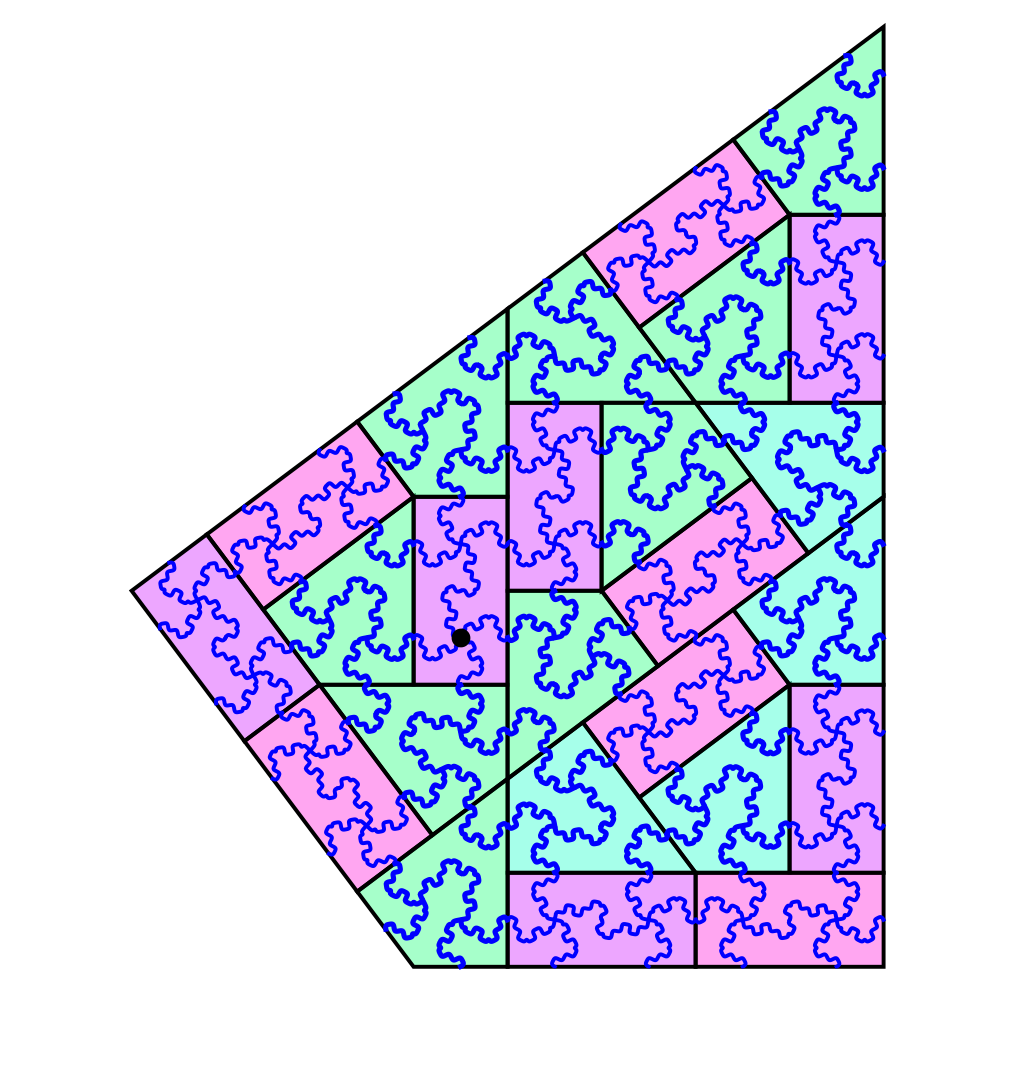}\hspace{.5cm} \includegraphics[width=3.35in]{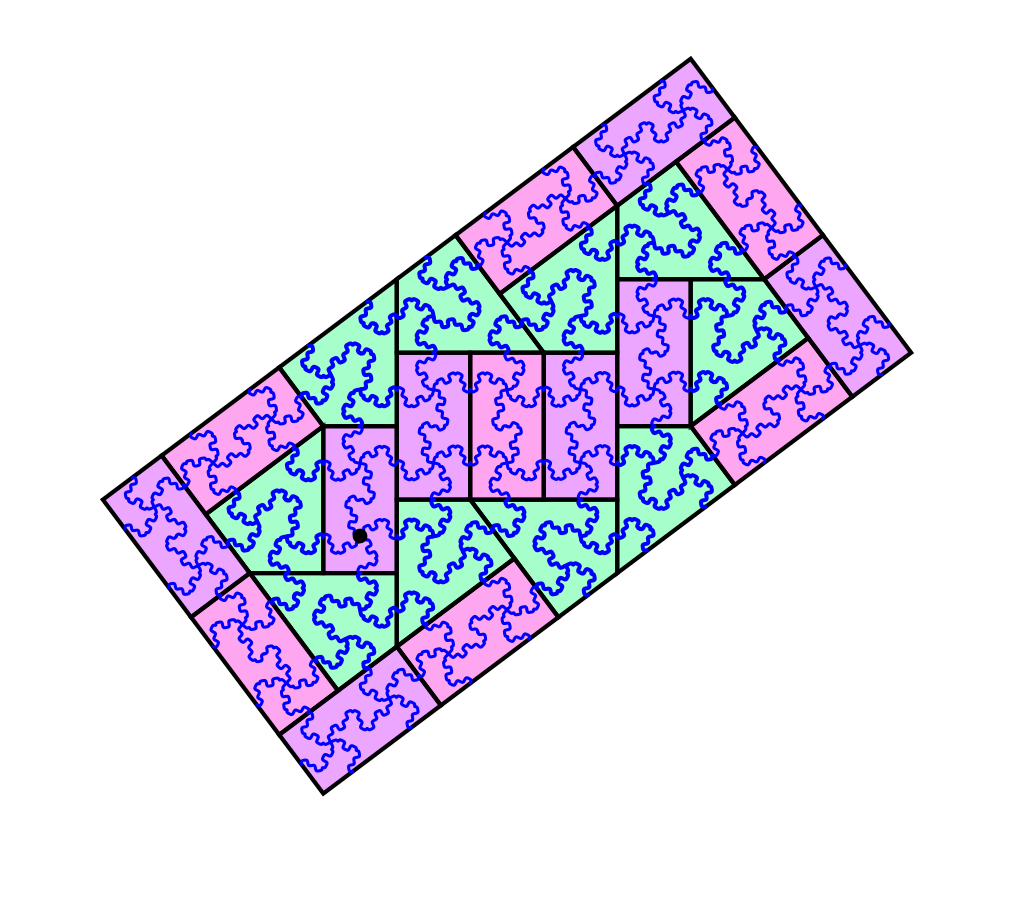}
\caption{Marking the level-1 kite and domino.}
\label{marked.patch}
\end{figure}

Figure \ref{marked.patch} shows the result of marking the kites and dominos this way in substituted kite and domino tiles. 
The fractal marking of any kite or domino will join with the marking of its neighbors at the side control points forming a fractal connection between their central control points.    The fractal connections encircle closed regions, and when such a region has no fractal in its interior we call it a  {\em pinwheel fractile}.  Working by hand, we were not sure how many fractiles to expect and feared there could be hundreds.  We wrote computer code that generated further iterates of the marked kite-domino substitution and counted the fractiles that appeared in the images. By doing this we were relieved to find that there are 13 tile types up to reflection, and 18 tile types when reflection is considered distinct. (We defer temporarily the proof that we have exhausted all possibilities).  Of the 13 tile types, 10 are visible in Figure \ref{marked.patch}; the remaining three types appear after  one more kite-domino substitution.   

\begin{figure}[ht]
\raisebox{.6in}{\includegraphics[width=.9in]{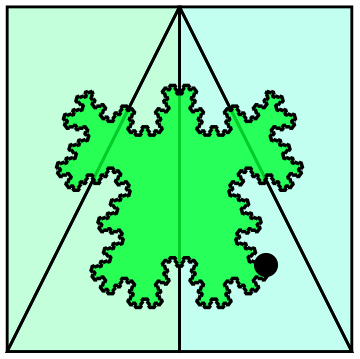}}\hspace{.5cm} \includegraphics[width=1.35in]{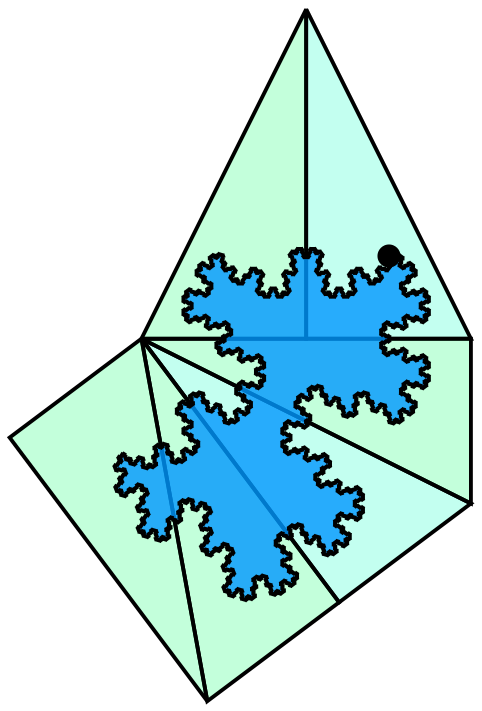}
\hspace{.5cm} \raisebox{.1in}{\includegraphics[width=1.5in]{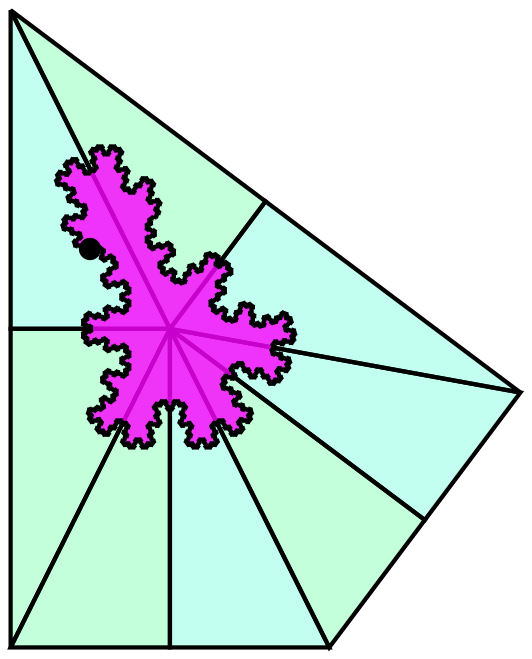}}
\caption{Three representative fractiles, in standard position, as they arise in their triangle patches.}
\label{representative.fractile.patches}
\end{figure}

Each fractile arises inside a patch of pinwheel triangles that is, except in a few cases, unique.  Knowing these patches is essential for writing computer code to generate the images, and for figuring out how to inflate and subdivide the fractiles.
In Figure \ref{representative.fractile.patches} we show three representative fractiles as they arise in their pinwheel triangle patches.
For readers who wish to get their hands dirty by drawing pinwheels of their own we have included a triangle in standard position in each patch and have marked the origin with a dot.   When a choice of which triangle to standardize needed to be made, we did so based on convenience for our computer code.  

\begin{figure}[h]
\includegraphics[width=4in]{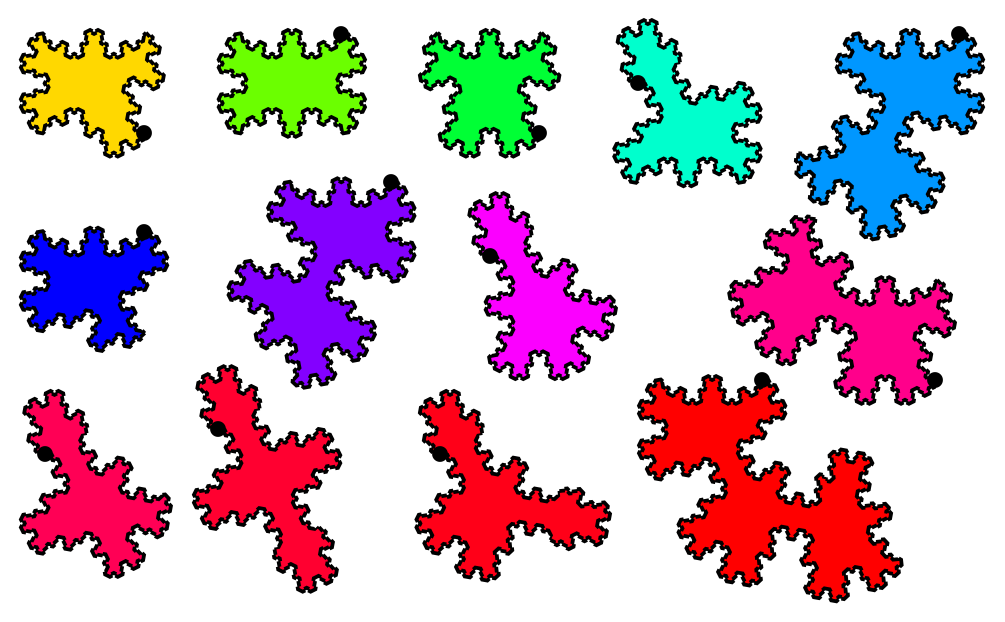}
\caption{The fractiles. For each tile the control point of the standard triangle is marked.}
\label{fractiles}
\end{figure}

 All thirteen fractiles are shown in Figure \ref{fractiles}, in order of relative frequency with the most frequently seen tile first.  (We take the frequency of the tiles to mean the average number of times the tile appears in any orientation, per unit area.) 
Using the control points as a reference the reader can see that they all fall into exactly one of the classes shown in Figure \ref{representative.fractile.patches}.   

\subsection*{Equivalence of tiling spaces}
\label{equiv.tiling.spaces}

We can mark any pinwheel tiling  $T$ in $X_P$ with the kite-domino method,
producing a new tiling made of pinwheel fractiles.  It is clear that doing this for every tiling in $X_P$ will produce a translation-invariant set of tilings, which in turn forms a tiling space that we denote $X_F$.  
Since the central control points of pinwheel triangles are exactly the locations where the continuations meet the aorta, the vertex set of a fractal pinwheel tiling and the set of central control points of the corresponding pinwheel triangle tiling coincide.

The pinwheel tiling space $X_P$, the kite-domino tiling space $X_{KD}$, and the fractal pinwheel tiling space $X_F$ are all mutually locally derivable.   The equivalence of the first two is in \cite{BFG1}; to complete the assertion we show that $X_{KD}$ and $X_F$ are mutually locally derivable.   
The fractile boundaries in any pinwheel fractal tiling  in $X_{F}$ are locally identifiable as kite or domino markings since their vertices are control points: if the vertex is degree 3, it is inside a kite; if it is degree 4, it is inside a domino.  Thus every tiling in $X_F$ locally determines a tiling in $X_{KD}$ and so $X_{KD}$ is locally derivable from $X_{F}$.   Conversely, any tiling in $X_{KD}$ can be marked as in Figure \ref{aortas.ventricles.kd}.  Once this is complete, kite-domino patches of radius 5 or smaller determine which fractile covers any given point in $\R^2$, since the largest fractile comes from a patch of kites and dominoes that has a diameter less than 5.     This means that $X_{F}$ is locally derivable  from $X_{KD}$ and completes the proof that all three tiling spaces are mutually locally derivable.  
\section*{The fractile substitution}
\label{basic.proofs}

\begin{figure}[ht]
\raisebox{.4in}{\includegraphics[width=.75in]{ghostpatch} }
\raisebox{.75in}{\text{\LARGE $\to$}}
\includegraphics[width=2.3in]{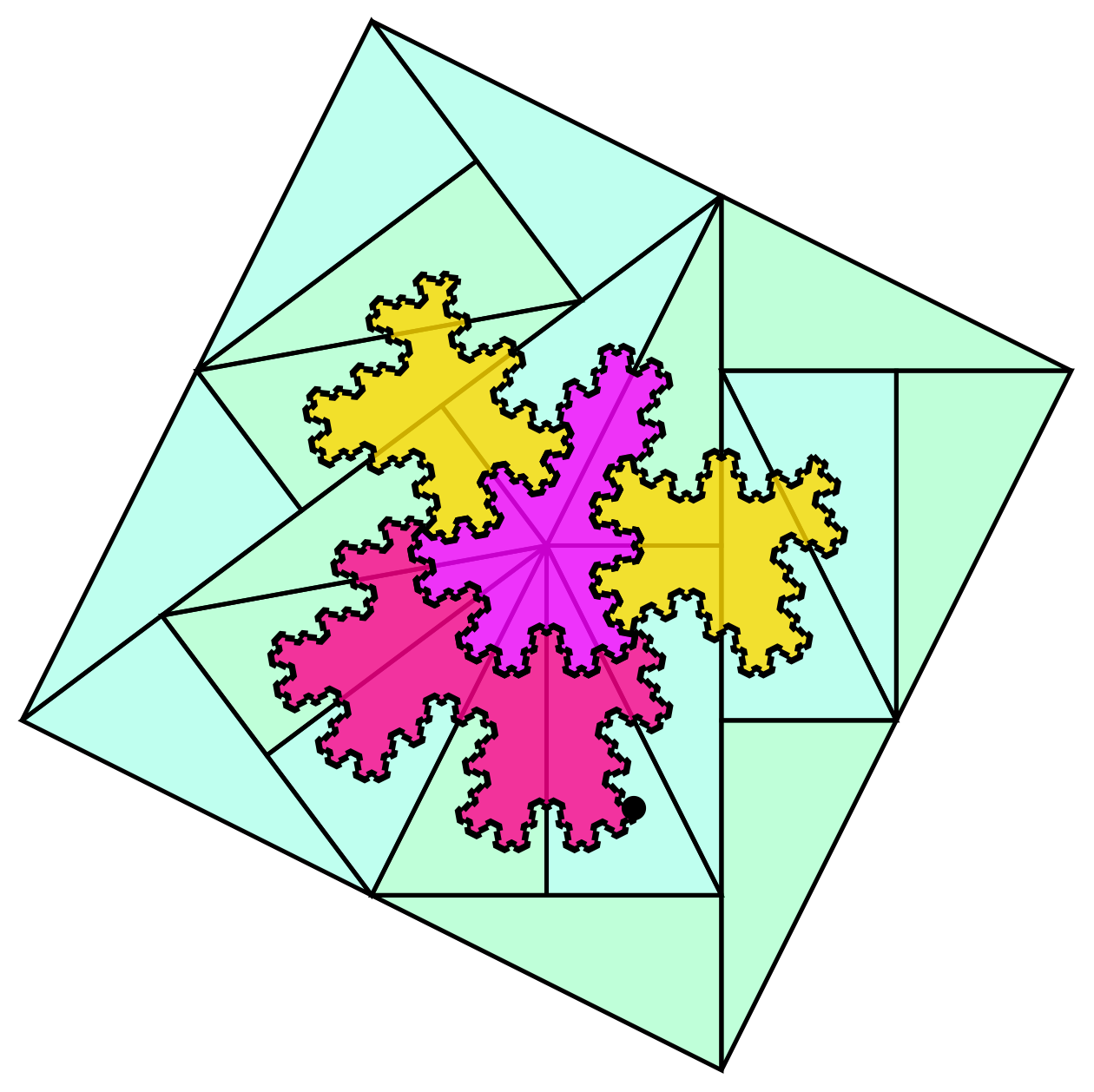}
\caption{How the ``ghost" fractile inherits its substitution rule.}
\label{inflate.ghost}
\end{figure}

The fact that each pinwheel fractile arises from a finite patch of pinwheel triangles means that each pinwheel fractile inherits a substitution from the pinwheel tiles that created it.   (It also inherits an equivalent one from the kite-domino substitution.) For a few of the fractile types, the triangle patch that creates it is not unique because the fractal markings at all four right angles of the domino tile create congruent regions (see Figure \ref{mark.domino.ventricle}).  However, it is easy to check that the pinwheel substitution induces congruent markings on the interiors of these regions, which implies that the substitution induced by the pinwheel substitution on the fractiles is well-defined.    In Figure \ref{inflate.ghost} we demonstrate how the substitution is induced on the fractile we call the ``ghost" (note that we include in this image only the portion of the kite-domino markings that lie inside the inflation of the ghost's boundary).    

\begin{figure}[ht]
\includegraphics[width=4.5in]{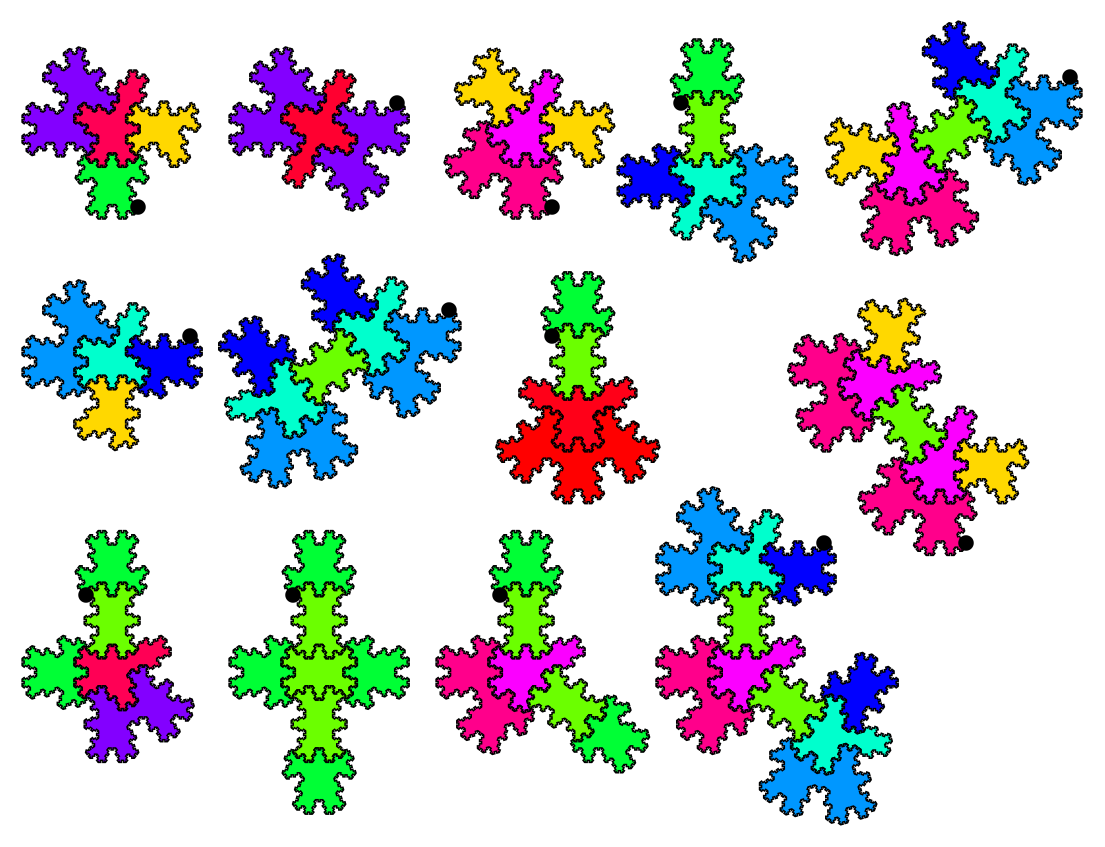}
\caption{Substituting the 13 fractile types.}
\label{fractile.subs}
\end{figure}

In Figure \ref{fractile.subs}, we show the substitutions of all thirteen basic pinwheel fractiles.    We would like to emphasize the remarkable fact that the boundaries of the tiles shown in Figures \ref{fractiles} and \ref{fractile.subs} are perfectly scaled versions of one another.    No additional detail is gained or lost because the tile boundaries are built from the aorta, and the aorta is a true fractal.

We can now argue that the list of fractiles we show in Figure \ref{fractiles} is complete.  Since the substitution rule shown in Figure \ref{fractile.subs} is self-contained, it defines a translation-invariant substitution tiling space $X_F'$ in the same way that the original pinwheel substitution generated the tiling space $X_P$.   $X_F'$ is measure-theoretically the same space as the space of pinwheel fractile tilings $X_F$ defined above; it is actually possible to show that they are exactly the same space.  For, $X_F'$ is mutually locally derivable from a subspace of $X_P$ which must be translation-invariant since $X_F'$ is.  Since $X_P$ is uniquely ergodic the subspace corresponding to $X_F'$ must have measure 0 or 1.  Since it is not of measure 0, it must correspond to a measure-1 subset of $X_P$. 

\subsection*{Fixed, periodic, and symmetric points in $X_F$}
\begin{figure}[h]
\raisebox{1.6in}{\includegraphics[width=.75in]{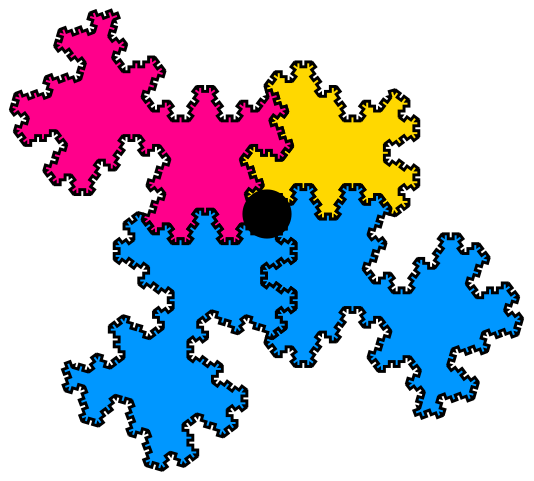} }
\raisebox{1.9in}{\text{\LARGE $\to$}}
\raisebox{1.1in}{\includegraphics[width=1.55in]{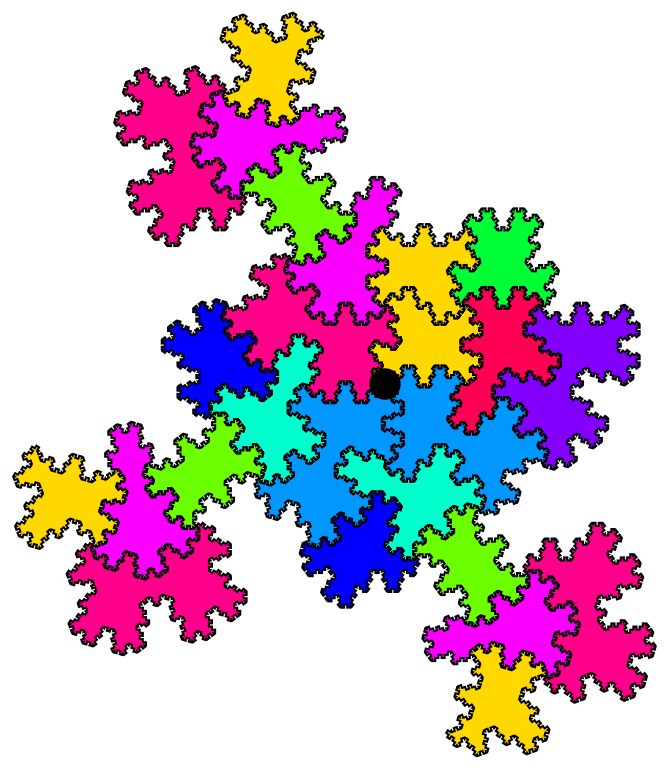}}
\raisebox{1.9in}{\text{\LARGE $\to$}}
\includegraphics[width=3.3in]{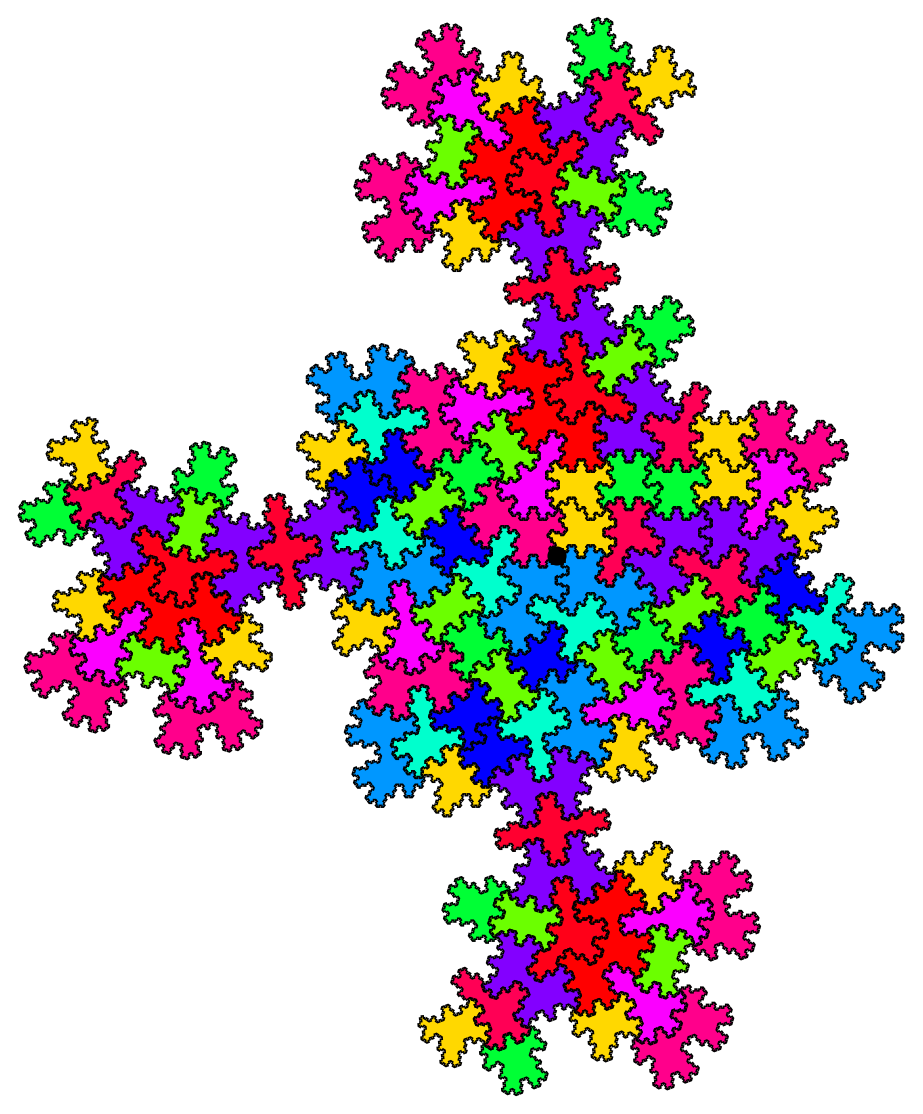}
\caption{Generating a fixed point of the fractile substitution.}
\label{origin.fixed}
\end{figure}

The substitution-invariance of the standard triangle (see Figure \ref{orig.pin.subs}) means that the fractile substitution also admits a fixed point (i.e., a self-similar tiling) as shown in Figure \ref{origin.fixed}. Of course any rotation of the standard pinwheel triangle also leads to a self-similar tiling.   The reflection of the standard pinwheel triangle across the $y$-axis almost leads to a substitution-invariant tiling too, but not quite:  the substitution of the reflected standard tile has the reflected standard tile at the origin, but it is rotated clockwise by the pinwheel angle $\phi = \arctan(1/2)$.  

\begin{figure}[ht]
\raisebox{1.4in}{\includegraphics[width=.625in]{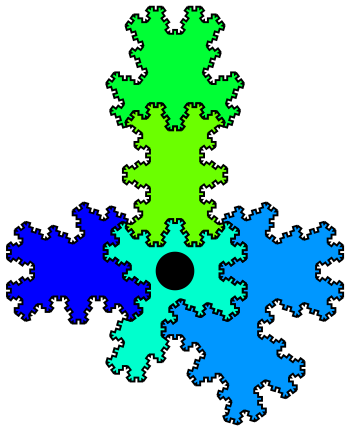} }
\raisebox{1.65in}{\text{\LARGE $\to$}}
\raisebox{.9in}{\includegraphics[width=1.45in]{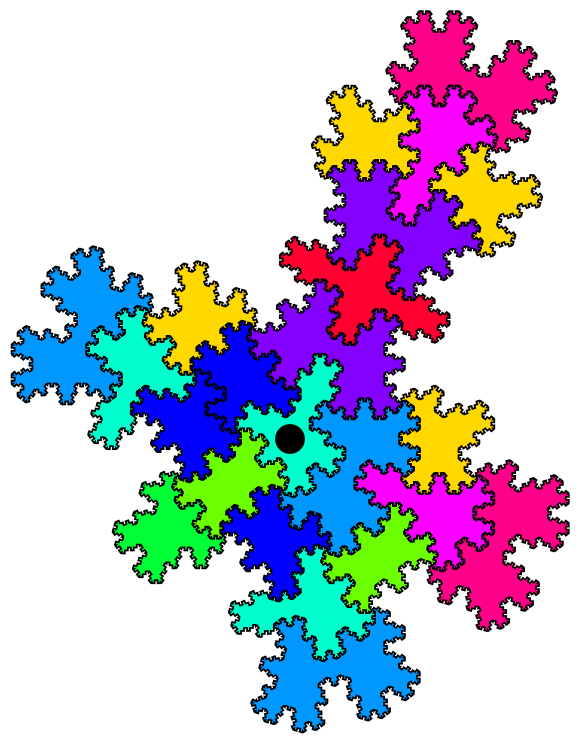}}
\raisebox{1.65in}{\text{\LARGE $\to$}}
\includegraphics[width=3.6in]{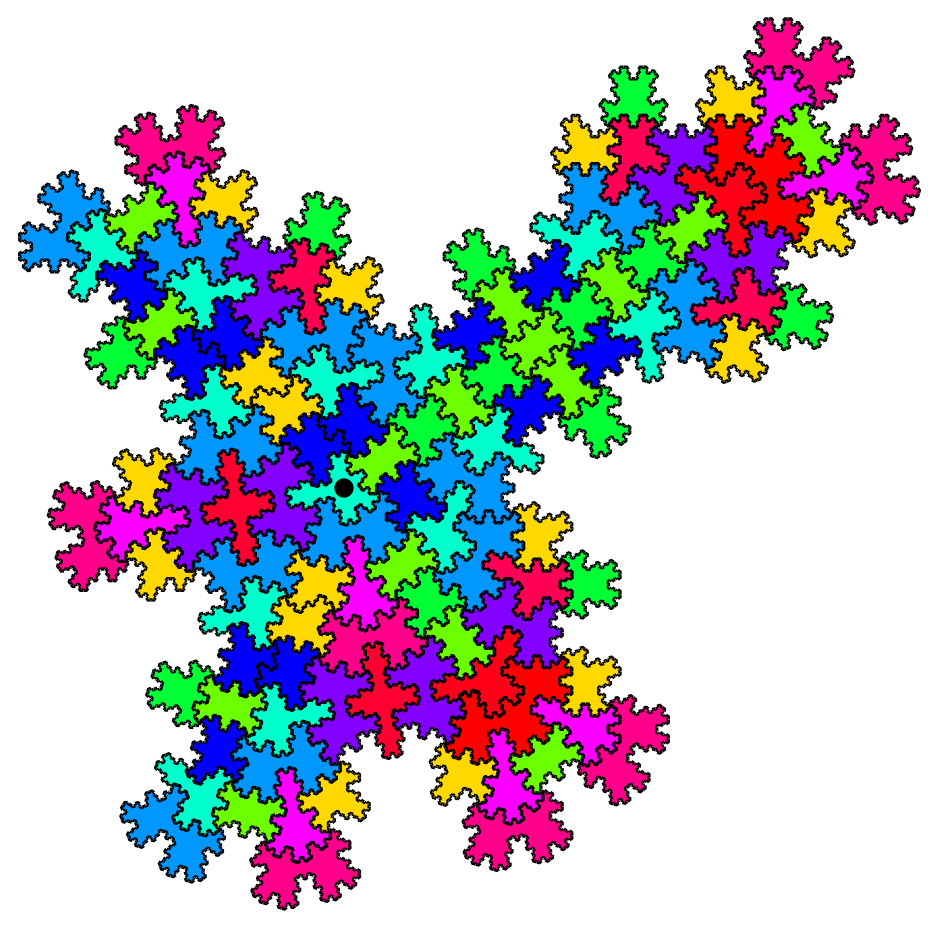}
\caption{Generating a tiling fixed by substitution and reflection.}
\label{origin.fixed.reflection}
\end{figure}

 One can check that there are no other fixed points by noticing that a tiling is invariant only if its patch at the origin is fixed under substitution.  This implies that the tiles in such a patch must, under substitution, contain themselves.  This happens only for the tiles pictured in Figure \ref{origin.fixed}.   It is interesting to note, however, that the fourth, sixth, and tenth prototiles contain reflections of themselves in their substitutions.  Thus the tilings they create are fixed when a combination of substitution and reflection are applied. Figure \ref{origin.fixed.reflection} develops how this looks with the fourth fractile at the origin.

\begin{figure}[ht]
\raisebox{1.65in}{\includegraphics[width=.75in]{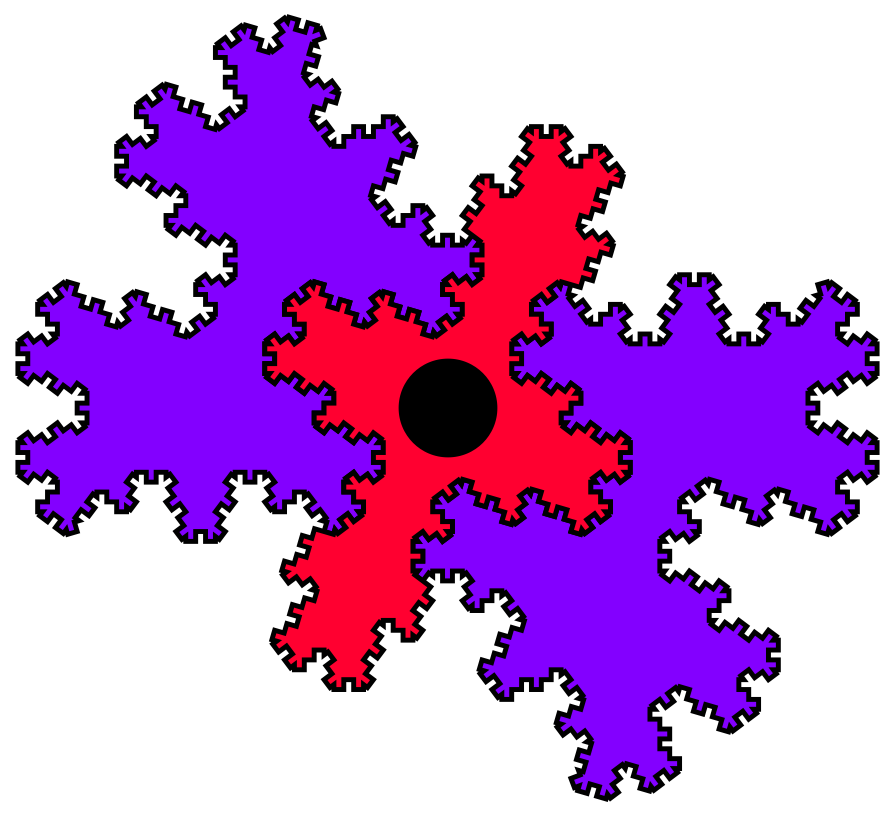} }
\raisebox{1.95in}{\text{\LARGE $\to$}}
\raisebox{1.1in}{\includegraphics[width=1.75in]{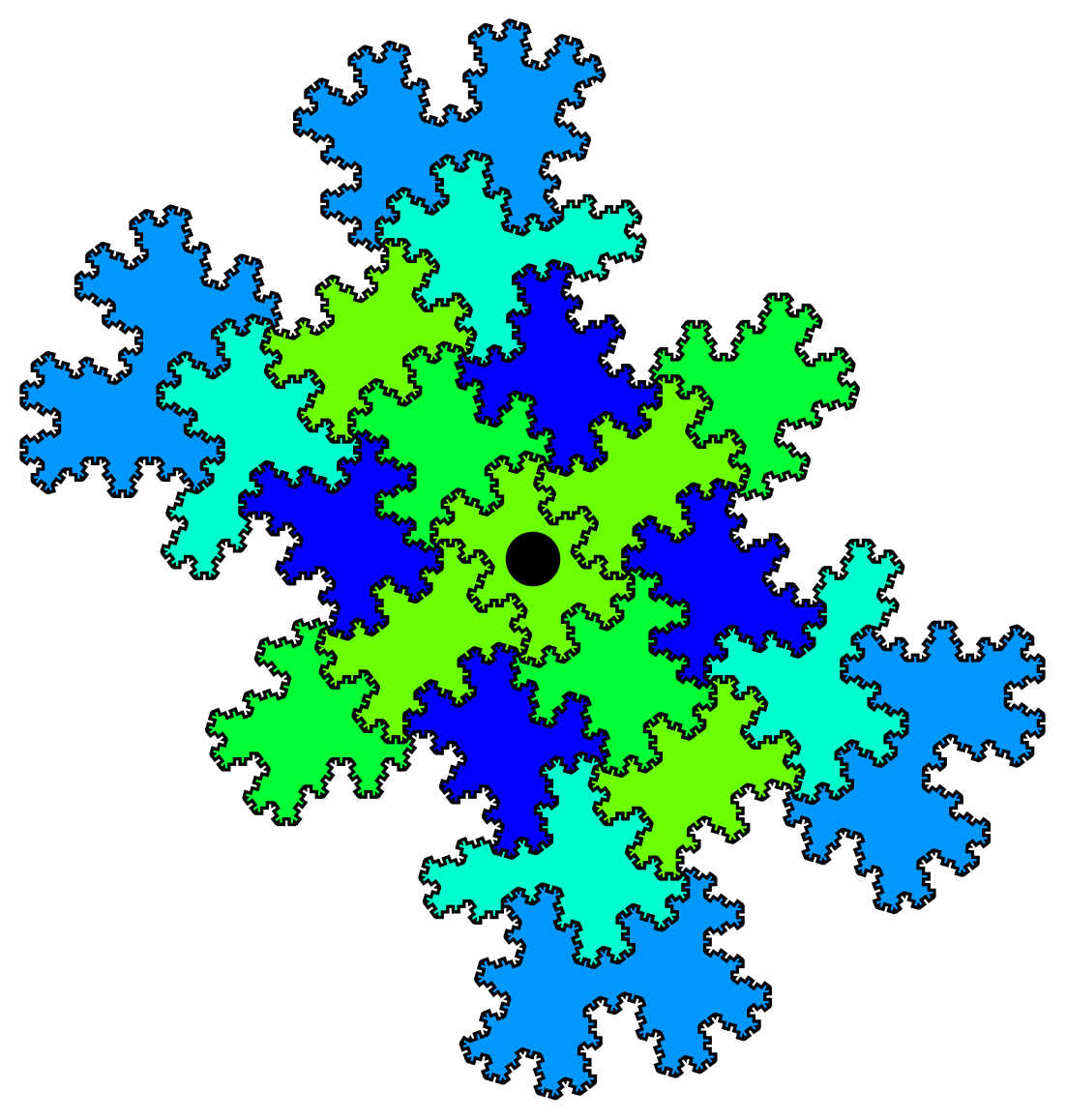}}
\raisebox{1.95in}{\text{\LARGE $\to$}}
\includegraphics[width=3.2
in]{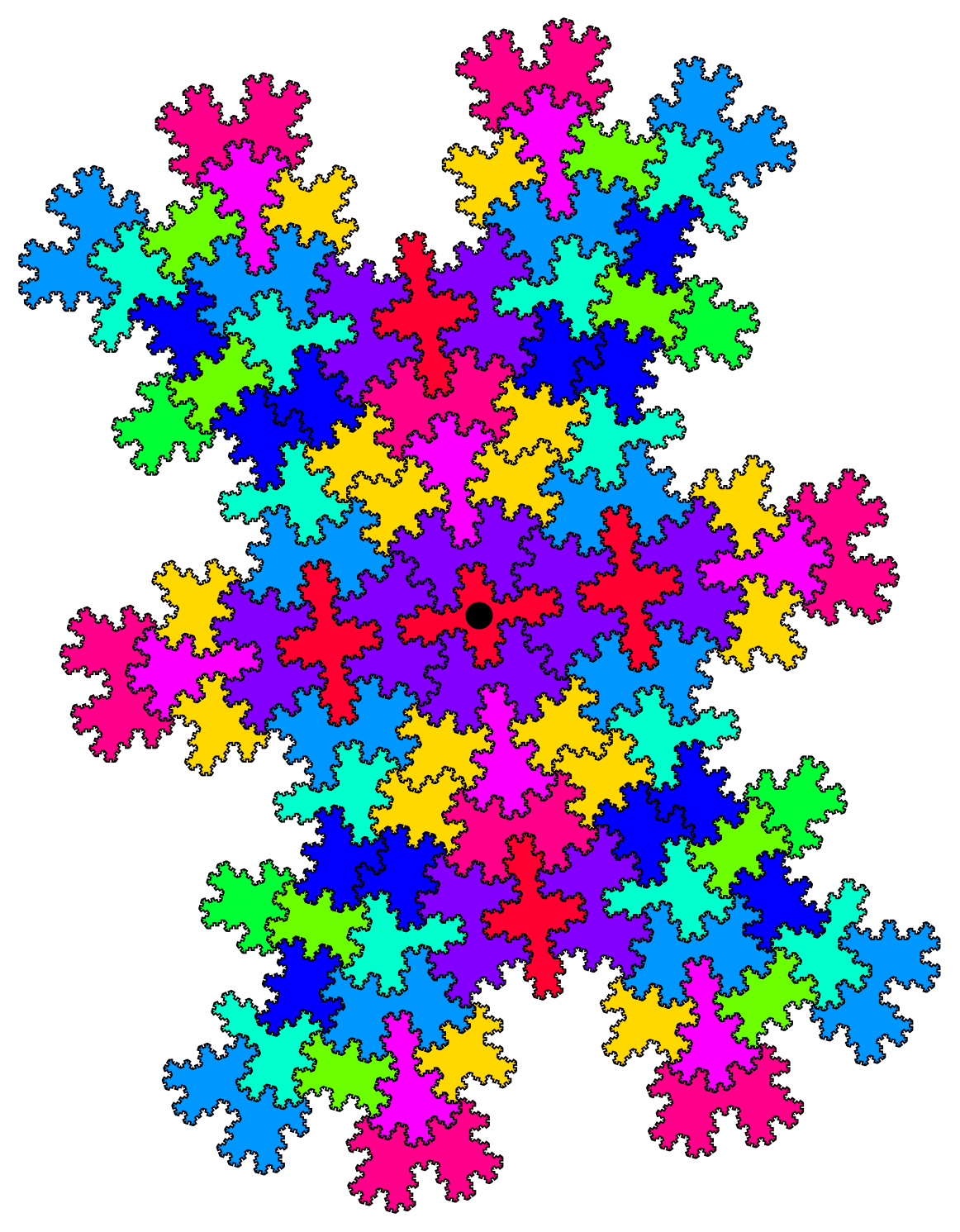}
\caption{Period 2 up to rotation by $\phi$; also invariant under rotation by $\pi$ and reflection.}
\label{origin.bat}
\end{figure}

There are six pinwheel triangle tilings that are fixed under rotation by $\pi$.
Two are invariant under reflection as well; these two are the images of each other under the original pinwheel substitution and are thus period-2 under substitution and rotation by 2$\phi$.   The corresponding tilings in $X_F$ have the second and the eleventh fractile types at the origin and are pictured in Figure \ref{origin.bat}.  

\begin{figure}[ht]
\raisebox{1.4in}{\includegraphics[width=.65in]{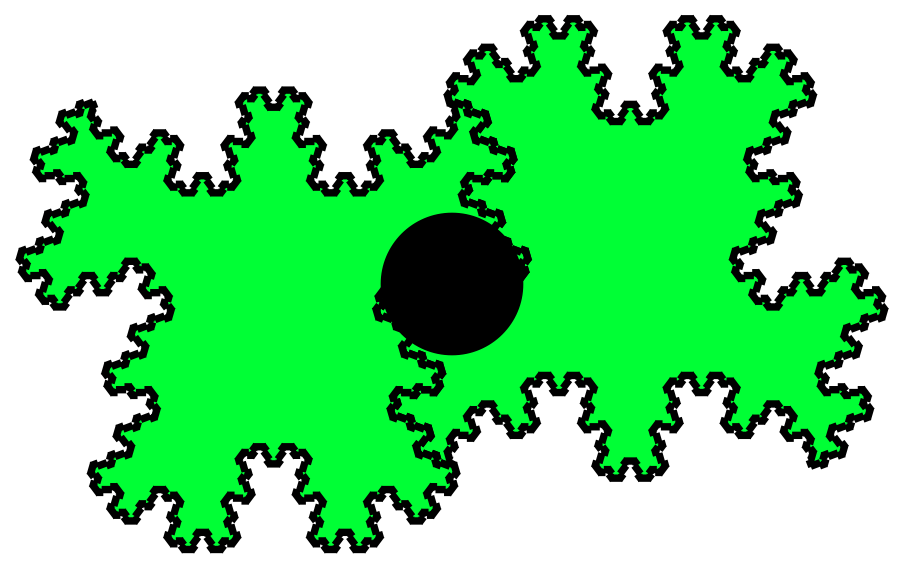} }
\raisebox{1.55in}{\text{\LARGE $\to$}}
\raisebox{1.0in}{\includegraphics[width=1.5in]{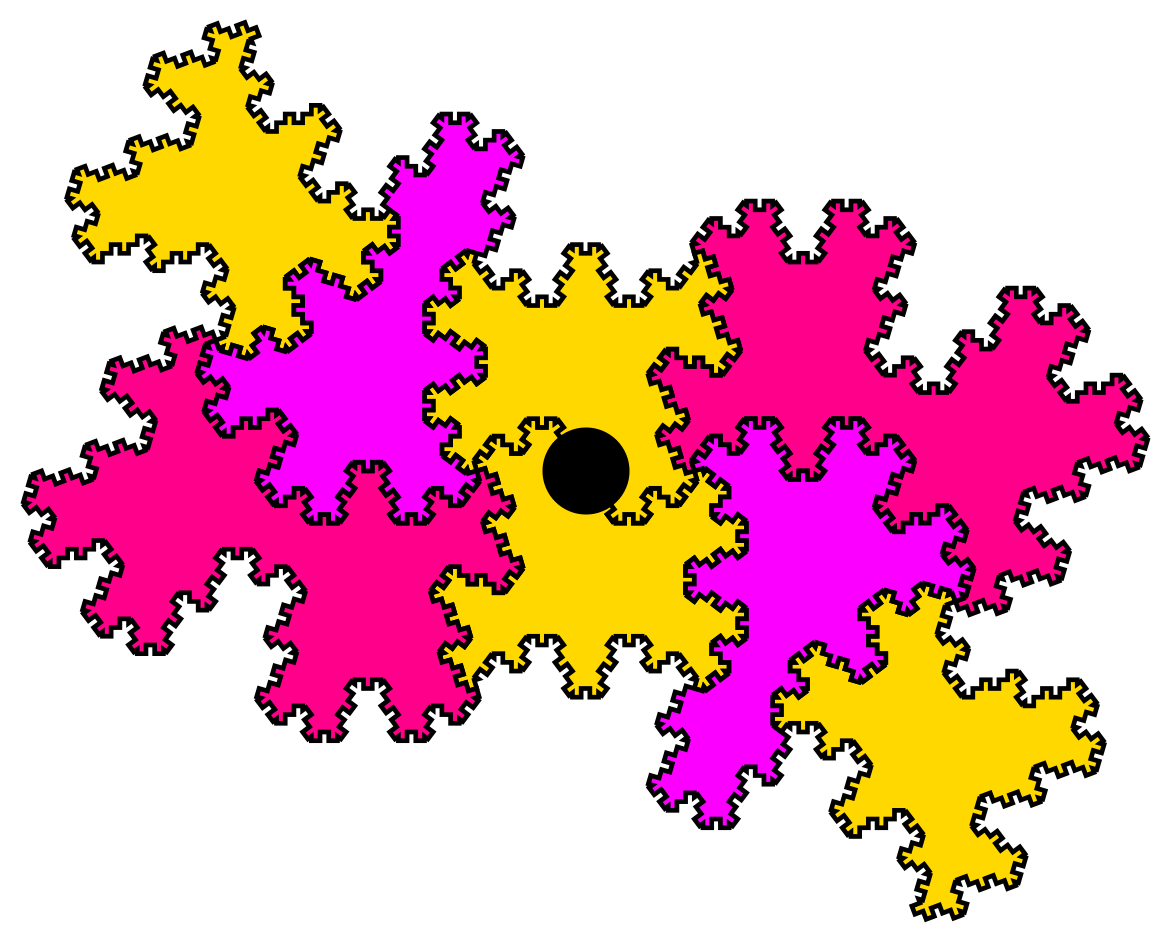}}
\raisebox{1.55in}{\text{\LARGE $\to$}}
\includegraphics[width=3in]{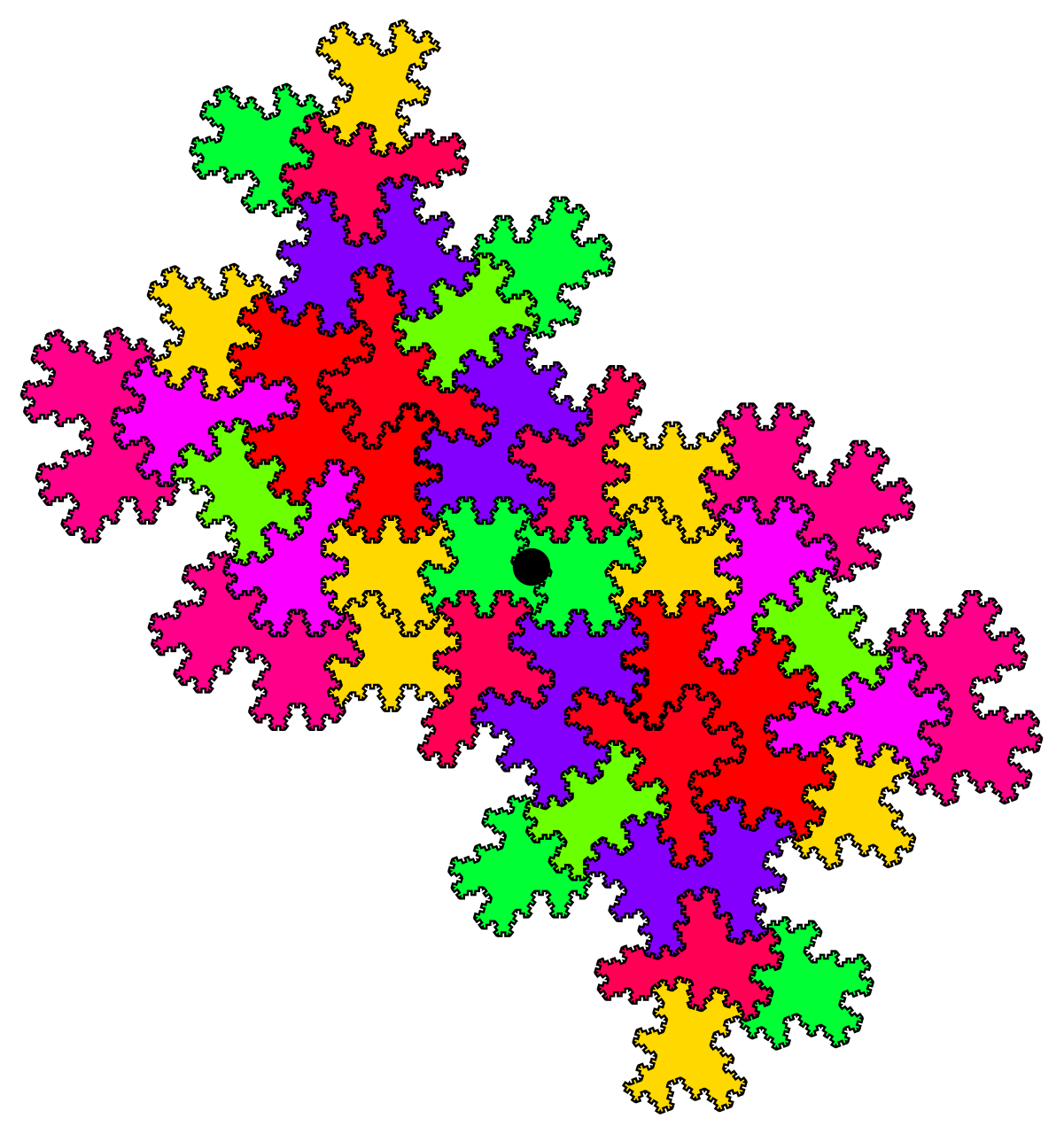}
\caption{Period 4, up to rotation by $\phi$; invariant under rotation by $\pi$.}
\label{origin.ghost}
\end{figure}

The other four pinwheel triangle tilings that are fixed under rotation have the center of a domino tile, or its image under substitution, at the origin.   These are not symmetric by reflection and make a period-4 sequence under substitution plus rotation, or a period-2 sequence under substitution plus reflection across an appropriate axis.   See Figure \ref{origin.ghost} for the beginning of the corresponding tilings in $X_F$.

\section*{Properties}
\label{advanced.proofs}

\subsection*{Basic properties of the fractiles}
A useful tool for analyzing substitution rules is the substitution matrix $A$ whose entry $A_{ij} $ is the number of tiles of type $j$ in the substitution of tile $i$ (cf. \cite{Sol.self.similar,Robinson.ams} for results used in this section).  Since $A$ is a nonnegative integer matrix, the largest eigenvalue $\lambda$ is real by the Perron-Frobenius theorem.  In fact, $\lambda$ is the area expansion of the substitution and its left eigenvector represents the relative areas of the tiles.   (In our case $\lambda=5$.) Moreover,  a properly scaled right eigenvector represents the relative frequency with which each tile appears, where relative frequency is the number of occurrences  per unit area.

We can choose whether or not to distinguish between tiles that are reflections of each other, giving us either 13 or 18 prototiles.   Although this affects the size of $A$, it does not affect the eigenvector analysis particularly much.   Obviously, the eigenvector representing the relative areas will give us the same relative sizes of tiles in each case.   And since reflections of tiles happen equally often, we find that when they are taken into account, the relative frequency is halved.   Thus, if we were to consider reflections to be distinct, the most frequently seen tile would no longer be the first tile shown in Figure \ref{fractiles} since it would appear half as often, and its reflection would also.    When we consider the number of prototiles to be 13, we compute the vector of relative frequencies to be approximately $(.1412, .1225, .1039, .1, .1, .1, .0843, .0784, .0784, .0353, .0245, .0157, .0157)$.
These numbers gave us the order we used to display the tiles in Figure \ref{fractiles}.

 At first it may seem surprising that the fractile areas are whole multiples of 1/5, but this fact can either be seen geometrically by looking at the kite and domino markings of Figure \ref{aortas.ventricles.kd} or by eigenvector analysis.   (The geometric argument begins by noticing that the aorta cuts its triangle exactly in half).  The areas of the tiles in the order shown in Figure \ref{fractiles} are: 1, 1, 1, 6/5, 9/5, 1, 9/5, 6/5, 9/5, 6/5, 7/5, 7/5, and 13/5. (Note that the area of a pinwheel triangle is also 1).  Since the aorta is the limit set of an IFS with a (linear) contraction factor $\sqrt{5}$ that uses three functions, its fractal dimension is $\ln{3}/\ln{\sqrt{5}} \approx 1.365$, thus the boundaries of the fractiles have this dimension also.  The tiles have rational area but irrational boundary dimension!

\subsection*{Rotational property}
\label{rot.props}

Since any pinwheel tiling features triangles in infinitely many different orientations, it is clear that any fractile must appear in any tiling in $X_F$ in infinitely many orientations also.  By the equivalence of the tiling spaces $X_P$ and $X_F$ we know that the pinwheel fractal tilings must also be ``statistically round"  in the same sense as the pinwheel triangle tilings.  All these orientations are wound up together in every copy of the aorta in an intriguing way.

 \begin{thm}
 \label{aorta.rotation}
 For every $N > 0$ there is a connected subset of the aorta, copies of which appear in at least $N$ distinct rotations inside the aorta.  Moreover, the set of all relative orientations that occur is uniformly distributed in $[0,2\pi]$.
 \end{thm}
 
 \begin{proof}
We will show that for each $N>0$ there is an $n\in \N$ for which the level-$n$ pinwheel supertile in standard position contains triangles that intersect the aorta and are in at least $N$ different relative orientations.  Applying the matrix $M_P^{-n}$ to this supertile will take the aortas of these triangles to the desired connected subsets of the aorta of the triangle in standard position.

We refer to Figure \ref{aorta.markings} for the level-1, 2, and 3 pinwheel supertiles in standard position.  In particular, notice that in the level-2 supertile, the triangle in standard position shares the vertex $(-.5, 1.5)$ with a triangle $t_1$ that is its rotation by $2\phi = 2\arctan(1/2)$ clockwise around that vertex.  After two more iterations of the substitution, we will have those two triangles, plus the two triangles in the same location in the substitution of $t_1$.   Of those, the image of the standard triangle is at the same orientation as $t_1$, but the image of $t_1$ is a rotation by another $2\phi$ clockwise.   So in this level-$4$ supertile, we have triangles at rotations of $0,2 \phi, $ and $4\phi$, and these tiles lie on the aorta.   In the level-$6$ supertile, we gain another triangle along the aorta, providing 4 distinct orientations.   In this way we see that if we need $N$ orientations, we must pass to a level-$(2N-2)$ supertile.

The second part of the theorem follows since $2\phi$ is irrational with respect to $\pi$, and thus the set $\{m (2\phi) \text{ such that} \, \, \, m \in \N\}$ is uniformly distributed mod $2\pi$.   
 \end{proof}

\begin{figure}[ht]
\includegraphics[width=4in]{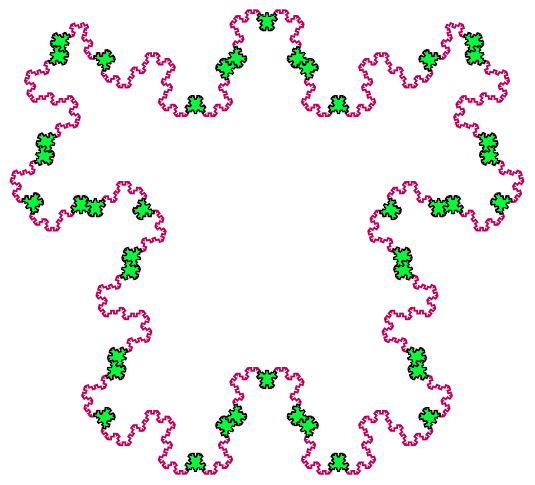}
\caption{Ghost tiles on the boundary of any level-4 ghost supertile}
\label{ghostdangle}
\end{figure}

We conclude with a fun side effect of the rotational and border-forcing properties.  Every fractal pinwheel tiling can be decomposed into level-$N$ supertiles for any $N$, and when $N$ is large so is the number of (level-0) fractiles the intersect the level-$N$ boundaries.  Our rotational theorem can be interpreted as:  if a certain fractile type (for instance, the ghost fractile) intersects the boundary of level-$N$ supertiles, then as $N$ increases to infinity, that fractile will dangle off the supertile boundaries in an unbounded number of orientations.    Since with very little modification we have a border-forcing substitution, we know that the way that these tiles dangle off will be identical every time a particular level-$N$ supertile appears.    We leave you with Figure \ref{ghostdangle}, which shows all of the ghost fractiles that intersect the boundary of any level-$4$ ghost supertile in any fractal pinwheel tiling, almost as children to a larger parent.  Although $N=4$ is not very large, we begin to see the many different angles in which these offspring appear.

 \end{document}